\documentclass[11pt,oneside,a4paper,final]{article}
\usepackage{amsmath,amssymb,makeidx,amsfonts,latexsym,enumerate,amsthm}
\usepackage{cite}

\newtheorem{lemma}{Lemma}[section]
\newtheorem{theorem}{Theorem}[section]

\begin{document}

\author{Karl-Olof Lindahl\\
School of Computer Science, Physics and Mathematics\\
Linnaeus University, 351 95, V\"{a}xj\"{o}, Sweden\\
\texttt{karl-olof.lindahl@lnu.se}\\
\and Michael Zieve\\
Department of Mathematics\\
University of Michigan, 530 Church Street\\
Ann Arbor, MI 48109-1043, USA\\
\texttt{zieve@umich.edu}
}

\title{On hyperbolic fixed points in ultrametric dynamics\footnote{\textbf{Published in \emph{p-Adic Numbers, Ultrametric Analysis and Applications}}, Vol. 2, No 3, pp. 232--240, 2010}}

\date{July 8, 2010}

\maketitle

\begin{abstract}
Let $K$ be a complete ultrametric field. We give lower and upper bounds for the size of linearization discs for power series over $K$ near hyperbolic fixed points. These estimates are maximal in the sense that there exist examples where these estimates give the exact size of the corresponding linearization disc.  In particular, at repelling fixed points, the linearization disc is equal to the maximal disc on which the power series is injective.

\end{abstract}

\vspace{1.5ex}\noindent {\bf Mathematics Subject Classification
(2000):} 32P05, 32H50, 37F50

\vspace{1.5ex}\noindent {\bf Key words:} dynamical system,
linearization, conjugacy, ultrametric field

\section{Introduction}

In recent years there has been an increasing interest in ultrametric (non-Archimedean) dynamics, i.e. the study of iteration of rational maps over a complete ultrametric field $K$ or,  more generally, over the corresponding projective line $\mathbb{P}(K)=K\cup\{\infty \}$. Ultrametric dynamics can be viewed both as an analogue of complex dynamics, and as a part of algebraic dynamics, see e.g. the recent monographs \cite{AnashinKhrennikov:2009,Silverman:2007} and references therein.
An important issue in this study are the periodic points of such maps. As in complex dynamics, the character of a periodic point is determined by the modulus of the multiplier (the derivative at the periodic point). 
A periodic point is attracting, indifferent, or repelling if the modulus of the multiplier is less than one, equal to one, or greater than one respectively. In this paper we consider local linearization near attracting or repelling periodic points of period one. As we study local properties it is natural to identify the map with its taylor series about the fixed point. 
Recall that a locally convergent power series $f\in K[[x]]$ of the form
\[
f(x)=\lambda x+ \text{higher order terms},
\]
is said to be locally linearizable at the fixed point at the origin if there is a convergent power series $g$ such that the semi-conjugacy (Schr\"{o}der functional equation)
\[
g\circ f(x)=\lambda g(x),
\]
holds on some non-empty disc about the origin. Such a semi-conjugacy always exists in the hyperbolic as well as in the one-dimensional $p$-adic non-parabolic ($\lambda$ not a root of unity) indifferent case as shown by Herman and Yoccoz \cite{HermanYoccoz:1981}.  In this paper we estimate the maximal disc on which a semi-conjugacy holds in the attracting and repelling case respectively. Furthermore, we also estimate the maximal disc on which the full conjugacy
 \[
 g\circ f\circ g^{-1}(x)=\lambda x,
 \]
is valid. Note that the conjugacy relation implies that for all iterates of $f$ we have 
$g\circ f^{n}\circ g^{-1}(x)=\lambda^n x$. As manifest in this paper the semi-conjugacy equation can hold on a strictly larger domain than that of the full conjugacy. The reason for this is that the full conjugacy does not allow $f$ to have roots of iterates, nor periodic points, whereas the semi-conjugacy might.
For example, we prove that for polynomials of the form 
$f(x)=\lambda x + a_2x^2$, $0<|\lambda |<1$, the semi-conjugacy holds on the disc
 $\left\{ |x|<1/|a_2|\right\}$, whereas the full conjugacy holds only on the strictly smaller disc 
 $\left\{ |x|<|\lambda|/|a_2|\right\}$. In this case the closest zero of $f$ to the origin sits on the sphere 
 $\left\{|x|=|\lambda|/|a_2|\right \}$, breaking the conjugacy there. 
 Moreover,  consider
 $f(x)=\lambda x + x^2/\lambda^2+x^3/\lambda ^3+\dots$, 
 $0<|\lambda| <1$. In this case $f$ converges for $|x|<|\lambda |$, the semi-conjugacy holds for
$|x|<|\lambda|^2 $, whereas the conjugacy holds for $|x|<|\lambda |^3$.
 On the other hand, concerning the corresponding repelling case we prove that for $f(x)=\lambda x + x^2/\lambda^2+x^3/\lambda ^3+\dots$,  the domain of 
 semi-conjugacy and full conjugacy both coincide with the domain of convergence of $f$; the disc $\{|x|<|\lambda |\}$.


The attracting case is summarized in Theorem \ref{lambda<1}, and the corresponding  repelling case  is given by Theorem \ref{theorem-repelling}.  Our results
holds for general complete ultrametric $K$, i.e. $K$ could be either some $p$-adic fileld like $\mathbb{C}_p$, or  some extension of a function field of Laurent series $\mathbb{F}((T))$ in positive or zero characteristic.

Estimates of the semi-conjugacy were previously obtained in the attracting case when $K$ is a finite extensions of the $p$-adic numbers $\mathbb{Q}_p$ by Lubin \cite{Lubin:1994}, and in the complex $p$-adic case 
$K=\mathbb{C}_p$ by Rivera-Letelier \cite{Rivera-Letelier:2003thesis}. In these works the the conjugacy is constructed as an iterative logarithm whereas our approach uses the ansatz of a power series $g$ and estimates of its coefficients. Similarly, in the repelling case, B\'{e}zivin 
\cite{Bezivin:2001} studied the conjugacy in the repelling case over $K=\mathbb{C}_p$, and gave a proof of the existence of the conjugacy. This result was then used to prove a result on the properties of the corresponding Julia set as mentioned below. Moreover, similar results to ours, concerning the convergence of the conjugacy, were stated by Zieve in the thesis \cite{Zieve:1996}. See Remark 1 and 7 below.

In the non-parabolic indifferent case there do not always exist a linearization over fields of positive charactersitic. (Paracolic germs are in general not linearizable over any ground field.) The indifferent case was studied in the equal characteristic cases \cite{Lindahl:2004,Lindahl:2009eq, Lindahl:2010} and in the $p$-adic case \cite{Ben-Menahem:1988,ThiranVerstegenWeyers:1989,
Arrowsmith/Vivaldi:1994,PettigrewRobertsVivaldi:2001,Khrennikov:2001a, Zieve:1996} later generalized in \cite{Lindahl:2009pp}, and by Viegue in the multi-dimensional 
$p$-adic case \cite{Viegue:2007}.


The paper is organized in the following way. First we recall some facts in section \ref{background} about the mappings that we study. In the following sections we state and prove the results in the  attracting and repelling cases respectively. The stated results are followed by some remarks concerning implications and related works.

Before closing this section we briefly discuss some general facts about the corresponding hyperbolic dynamics. As in the complex field case a rational map $R\in K(x)$ induces a partition of the state space $\mathbb{P}(K)$ into the Fatou set $F(R)$, the domain on which $R$ is equicontinous and loosely speaking the `region of order', and its complement, the Julia set $J(R)$ which is the `region of chaos', where small errors become huge after many iterations of $R$. There is a vast number of papers studying the properties of these sets, see e.g. \cite{Benedetto:2001a,Benedetto:2001b,Benedetto:2002, Benedetto:2006, Bezivin:2004b,Bezivin:2004a, Bezivin:2005,Hsia:2000, Rivera-Letelier:2003thesis, Rivera-Letelier:2000, FavreRivera-Letelier:2010}. In particular, attracting periodic points are all in $F(R)$, while repelling periodic points all belong to the Julia set $J(R)$. In fact, as proven by B\'{e}zivin \cite{Bezivin:2001}, the Julia set is the closure of the repelling periodic points. There are several characterizations of the Fatou set in terms of components. Rivera-Letelier \cite{Rivera-Letelier:2003thesis} showed that  components containing an attracting fixed point $x_0$, and all the elements that are attracted to $x_0$, may have a complicated `cantor-like' structure. Another very interesting fact is that, contrary to the complex field case, the Fatou set may contain
so called `wandering domains', i.e. there may exist components that are not preperiodic, as proven by Benedetto
\cite{Benedetto:2002wandering,Benedetto:2006}. Beside those works already mentioned, ultrametric hyperbolic dynamics has been studied in many other papers, see e.g. 
\cite{KingsberyLevinPreygelSilva:2009,Khrennikov:2001a,KhrennikovNilsson:2001, KhrennikovNilssonMainetti:2002,
DeSmedtKhrennikov:1997,
Khrennikov/Nilsson:2004,
Khrennikov:2003ryssbok,DragovichKhrennikovMihajlovic:2007,
KhrennikovMukhamedovMendes:2007,KhrennikovSvensson:2007,
Khrennikov:2003nauk, Svensson:2005,NilssonNyqvist:2004,Li:2002a} and references therein.

\section{Background}\label{background}

In this section we introduce the main definitions and recall some basic facts from ultrametric calculus. For a general treatment of ultrametric power series we refer to
\cite{Schikhof:1984,Escassut:1995,Robert:2000}.

Throughout this paper $p$ is a prime and the base field $K$ is a complete ultrametric field with residue field $\Bbbk$. This means that $K$ can be of any of the following three types;  the $p$-adic case in which the characteristic char $K=0$ and char $\Bbbk=p$ (the so called mixed characteristic case), or one of the equal characteristic cases in which char $K=$ char $\Bbbk$ is either zero or 
$p$. The main types in the equal characteristic case are the function fields described below.

\textbf{Example 1} (char $K=$ char $k=p$). Let $F_p$ be a field of characteristic $p>0$, e.g.\@ $F_p$ could be the finite field
$\mathbb{F}_{p^n}$ of $p^n$ elements,  an infinite field like the algebraic closure 
of $\mathbb{F}_{p^n}$, or the infinite field of rational functions over $\mathbb{F}_{p^n}$. 
Let $K=F_p((T))$ be the field of
all formal Laurent series in variable $T$, with coefficients
in the field $F_p$.  An element $x\in K$ is of the form
\begin{equation}\label{laurent series}
x=\sum_{i\geq i_0} x_iT^i, \quad x_{i_0}\neq 0, \text{ } x_i\in
\mathbb{F},
\end{equation}
for some integer $i_0\in \mathbb{Z}$. 
Given $0<\epsilon <1$, we define an absolute value $|\cdot |$ on $K$ such that $|T|=\epsilon$ and
\begin{equation}\label{definition absolute value}
\left |\sum_{i\geq i_0} x_iT^i\right |=\epsilon ^{i_0}.
\end{equation}
Hence, $\pi=T$ is a uniformizer of $K$. Furthermore, $K$ is complete with respect to $|\cdot |$ and, analogously to the $p$-adic numbers, can be viewed as the completion of the field of
rational functions $F_p(T)$ over $F_p$ with respect to the absolute
value defined by (\ref{definition absolute value}).
Note that $i_0$ is the order of the zero (or if negative, the order of the pole) of $x$ at $T=0$.  Note also that in this case the residue field $\Bbbk=F_p$.
 Moreover,  $|\cdot |$ is the
trivial absolute value on $F_p$, the subfield of $K$ consisting of all
constant series in $K$. As for the $p$-adic  numbers, we can construct a completion $\widehat{K}$ of an algebraic closure
of $K$ with respect to an extension of $|\cdot |$. 
Then $\widehat{K}$ is a complete, algebraically closed non-Archimedean field, and its
residue field $\hat{k}$ is an algebraic closure of $k=F_p$. It follows that
$\hat{k}$ has to be infinite. The value group $|\widehat{K}^*|$, that is the set of real numbers which
are actually  absolute values of non-zero elements of $\widehat{K}$, will consist
of all rational powers of $\epsilon$, rather than just integer powers of $\epsilon$ as in $|K^*|$.
In particular, the absolute value is discrete on $K$ but not on $\widehat{K}$.

\textbf{Example 2} (char $K=$ char $k=0$).
Let $F$ be a field of characteristic zero, e.g.\@ $F$ could be either $\mathbb{Q}$, $\mathbb{Q}_p$,
$\mathbb{R}$ or $\mathbb{C}$. Then define $K=F((T))$ and $|\cdot |$
as in the previous example. One major difference is that in this case the residue field $k=F$
must be infinite  since $F$ is of characteristic zero.

Given an element $x\in K$ and real number $r>0$ we denote by
$D_{r}(x)$ the open disc of radius $r$ about $x$, by
$\overline{D}_r(x)$ the closed disc, and by $S_{r}(x)$ the sphere
of radius $r$ about $x$. If $r\in|K^*|$ (that is if $r$ is
actually the absolute value of some nonzero element of $K$), we
say that $D_{r}(x)$ and $\overline{D}_r(x)$ are \emph{rational}.
Note that $S_r(0)$ is non-empty if and only if $r\in|K^*|$. If
$r\notin |K^*|$, then we will call $D_{r}(x)=\overline{D}_r(x)$ an
\emph{irrational} disc. In particular if $a\in K$ and $r=|a|^s$
for some rational number $s\in\mathbb{Q}$, then $D_{r}(x)$ and
$\overline{D}_r(x)$ are \emph{rational} considered as discs in
$\widehat{K}$.  Note that all discs are both open and closed as
topological sets, because of ultrametricity. However, power series distinguish between rational
open, rational closed, and irrational discs.

Let $K$ be a complete ultrametric field. Let $h$ be a power
series over $K$ of the form
\begin{equation*}
h(x)=\sum_{i=1}^{\infty}c_i(x-\alpha )^i, \quad c_i\in K.
\end{equation*}
Then $h$ converges on the open disc $D_{R_h}(\alpha )$ of radius
\begin{equation}\label{radius of convergence}
R_h =1/ \limsup |c_i| ^{1/i},
\end{equation}
and diverges outside the closed disc $\overline{D}_{R_h}(\alpha )$. The power series $h$ converges on the sphere $S_{R_h}(\alpha )$ if and only if
\[
\lim_{i\to\infty}|c_i| R_h ^i=0.
\]
For our purposes it is enough to consider the case $\alpha=0$ and $c_0=0$. Recall the following fact, which is a consequence of the Weierstrass Preparation Theorem.

\begin{lemma}\label{lemma one-to-one}

Let $K$ be algebraically closed and let
$h(x)=\sum_{k=1}^{\infty}c_kx^k$ be a power series over $K$.
\begin{enumerate}[1.]

 \item Suppose that $h$ converges on the rational closed disc
  $\overline{D}_R(0)$. Let $0<r\leq R$ and suppose that
  \[
   |c_k|r^k\leq |c_1|r\quad \text{ for all } k\geq 2 .
  \]
 Then $h$ maps the open disc $D_{r}(0)$ one-to-one onto
 $D_{|c_1|r}(0)$. Furthermore,
 \[
 d = \max\{k\geq 1:|c_k|{r}^k=|c_1| r\},
 \]
is attained and finite, and $h$ maps the  closed disc
$\overline{D}_{r}(0)$ onto $\overline{D}_{|c_1|r}(0)$ exactly
$d$-to-1 (counting multiplicity).

 \item Suppose that $h$ converges on the rational open disc
  $D_R(0)$ (but not necessarily on the sphere $S_R(0)$).
  Let $0<r\leq R$ and suppose that
  \[
   |c_k|r^k \leq |c_1|r\quad \text{ for all } k\geq 2 .
  \]
  Then $h$ maps $D_{r}(0)$ one-to-one
  onto $D_{|c_1|r}(0)$.
\end{enumerate}

\end{lemma}
 
A proof is given, in a more general setting, in \cite{Benedetto:2003a}.

Given the real number $\rho$, we will be interested in the following family of maps wtih multiplier $\lambda\neq 0$ fixing the origin in $K$.
\begin{equation}\label{definitionfamily}
\mathcal{F}_{\lambda,\rho}(K)=\left \{ f\in K[[x]]: f(x)=\lambda x+\sum_{i=2}^{\infty} a_i x^{i},  1/\sup_{i\geq 2} |a_i|^{1/(i-1)}=\rho
 \right \}.
\end{equation}
For future reference, we note that by definition 
\begin{equation}\label{boundf}
|a_i|\leq \left( 1/\rho\right)^{i-1}.
\end{equation}
To each such $f\in \mathcal{F}_{\lambda,\rho}(K)$ we will associate its radius of convergence
\begin{equation}
R_f=1/\limsup |a_i|^{1/i},
\end{equation}
and the maximal radius on which $f$ is one-to-one
\begin{equation}
\delta_f=\min\left\{ R_f, \inf_{i\geq 2}\left (|\lambda| /|a_i|\right )^{1/(i-1)}\right\} 
\end{equation}
respectively. 
In fact, by the lemma above
\begin{equation}\label{eq-disc-onetoone}
f: D_{\delta_f}(0)\to D_{|\lambda|\delta_f}(0)
\end{equation}
is a bijection. The radius of injectivity is discussed in more detail in the two sections below. 

\subsection{Injectivity -- Attracting case}\label{section injectivityattracting}
Let us now consider the attracting case $0<|\lambda|< 1$. Clearly $\delta_f\leq \rho\leq R_f$.
On the other hand, by (\ref{boundf}) we have $\delta_f\geq |\lambda|\rho$. 
Consequently,
\begin{equation}\label{boundfinjectivityattracting}
 |\lambda |\rho \leq \delta_f \leq \rho\leq R_f, \quad \text{if }  |\lambda| < 1.
\end{equation}
In particular, if $|a_i|=|\lambda |$ for all $i\geq 2$, then $\delta_f=\rho =R_f=1$. On the other hand, according to the example below,  we may also have $\delta_f=|\lambda| \rho$. This case is of special interest because, according to theorem \ref{lambda<1}, $|\lambda|\rho$ is a lower bound for the radius of the full conjugacy. 

\textbf{Example 3}.
Suppose that $0<|\lambda |<1$ and  $\rho= 1/|a_2|$, then 
\[
1/|a_2|\leq 1/|a_i|^{1/(i-1)} \text{ for all } i\geq 3.
\]
Moreover, $R_f\geq \rho$. Consequently,
\[
\delta_f=\inf_{i\geq 2}\left (|\lambda| /|a_i|\right )^{1/(i-1)}\leq 
|\lambda|/|a_2|=|\lambda|\rho.
\]
Recall that by (\ref{boundfinjectivityattracting}) we also have $\delta_f\geq |\lambda| \rho$.
It follows that the radius of injectivity
\begin{equation}\label{deltaequaltoradiusoflinearization}
\delta_f = |\lambda |\rho \quad \text{ if } \rho=1/|a_2|.
\end{equation}

\subsection{Injectivity -- Repelling case}\label{section injectivityrepelling}
We close the section on background material by considering the repelling case. 
In this case we have
\[
\inf_{i\geq 2}\left (|\lambda| /|a_i|\right )^{1/(i-1)}\geq \rho.
\]
On the other hand, for $|\lambda |>1$, we also have
\[
\inf_{i\geq 2}\left (|\lambda| /|a_i|\right )^{1/(i-1)}\leq |\lambda | \inf_{i\geq 2}\left (|\lambda|^{2-i} /|a_i|\right )^{1/(i-1)}\leq |\lambda |\rho,
\]
with equality if and only if $\rho=1/|a_2|$. Consequently,
\begin{equation}\label{boundfinjectivityrepelling}
\rho\leq \delta_f \leq \min\left \{R_f, |\lambda |\rho\right \} \quad\text{if  } |\lambda| > 1,
\end{equation}
and in particular, we may have $\delta_f=R_f$ in this case; e.g. if $|a_i|=|\lambda|$ for every 
$i\geq 2$, then $\delta_f=R_f=1=|\lambda |\rho$, since $\rho=1/|\lambda|$ in this case. 
That we may also have the other extreme in (\ref{boundfinjectivityrepelling}), $\delta_f=R_f=\rho$, follows from the following example.

\textbf{Example 4}.
Suppose that $|\lambda |>1$ and
\[
a_i=1/\lambda^{i} \quad \text{for } i\geq 2,  
\]
then

\[
 \limsup_{i \geq 1}|a_i|^{1/{i}}=\sup_{i \geq 2}|a_i|^{1/{i-1}}=1/|\lambda|, 
\]
and
\[
\inf_{i\geq 2}\left (|\lambda| /|a_i|\right )^{1/(i-1)}=
 \inf_{i\geq 2}\left (|\lambda|^{(i+1)/(i-1)}|\right )=|\lambda |.
\]
Accordingly,  $\delta_f=R_f=\rho= |\lambda|$ in this case.

\section{Attracting fixed points}\label{section-attracting}

In this section we will prove the following theorem.

\begin{theorem}\label{lambda<1}
Let  $K$ be a complete ultrametric field, and $f\in\mathcal{F}_{\lambda,\rho}(K)$, defined by 
(\ref{definitionfamily}), with 
$0<|\lambda|<1$. Then there is a
unique function $g$, defined on the open disc $\left\{|x|< \rho\right\}$ with
$g(0)=0$, $g'(0)=1$ such that the semi-conjugacy 
\begin{equation}\label{semi-conjugate-equation attracting}
 g\circ f(x)=\lambda g(x),
\end{equation}
is valid in the open disc $\left\{|x|< \rho\right \}$. The full conjugacy 
\[
g\circ f\circ g^{-1}(x)=\lambda x, 
\]
holds for $|x|<|\lambda|\rho$.
\end{theorem}
Before proving the theorem we make a few remarks. 

\noindent 
\textbf{Remark 1}.   
It follows from the proof below that one can easily obtain the estimate 
\[
 |b_k|\leq |\lambda|^{-k}\left (1/\rho\right)^{k-1},
\]
for the coefficients of the conjugacy. Consequently, as stated by Zieve \cite{Zieve:1996}, the semi-conjugacy  converges on the open disc of radius $\lambda \rho$ about the origin. In the proof below we improve this estimate by showing that in fact the coefficients of the conjugacy satisfy 
(\ref{boundbkattracting}) and converges on the open disc of radius $\rho$.

\noindent 
\textbf{Remark 2}.
A similar estimate, concerning the semi-conjugacy (\ref{semi-conjugate-equation attracting}) in Theorem \ref{lambda<1}, is known in the $p$-adic case
using the construction of an iterative logarithm 
$g=\lim_{n\to\infty}f^{\circ n}/\lambda^n$. 
Rivera-Letelier \cite{Rivera-Letelier:2000} (Proposition 3.3 p. 49) proved that, in the case 
$K=\mathbb{C}_p$, $g$ converges on the open disc of radius 
$r(f)=1/ \sup_{i\geq 1}|a_i|^{1/i} $ about the origin. This result was proven earlier by  Lubin \cite{Lubin:1994} in the special case that $K$ is a finite extension of $\mathbb{Q}_p$ and the coefficients of $f$ are all in the closed unit disc in $K$.

\noindent 
\textbf{Remark 3}. 
In view of (\ref{boundfinjectivityattracting}) the semi-conjugacy is valid on a disc of radius larger than or equal to that of the maximal disc on which $f$ is injective. In particular, if $|a_i|=|\lambda|$ for each $i\geq 2$, then as noted in section \ref{section injectivityattracting},
$\delta_f=\rho =R_f=1$. Hence, our estimate of the domain of semi-conjugacy is maximal in this case. Moreover, the maximal disc of semi-conjugacy is the same as that of the maximal disc about the origin on which $f$ is injective. In the next remark we conclude that there are cases in which the maximal disc on which $f$ is injective is strictly contained in the maximal disc of semi-conjugacy.

\noindent
\textbf{Remark 4}.
The estimate of the maximal disc of full conjugacy $\left\{|x|<|\lambda| \rho\right\}$ is sometimes the best possible; it may be that the radius of injectivity $\delta_f =|\lambda| \rho$, and consequently $f$ may have roots of iterates, i.e. $x$ such that $f^{n}(x)=0$, on the sphere $\left\{|x|=|\lambda |\rho\right\}$ breaking the conjugacy there; by the conjugacy relation 
$g\circ f^{n}\circ g^{-1}(x)=\lambda  ^nx$ and consequently $f$ could not have roots of iterates (nor periodic points) on the linearization disc. For example, in view of (\ref{deltaequaltoradiusoflinearization}), this is the case if $\rho =1/|a_2|$ as for quadratic polynomials.

\noindent
\textbf{Remark 5}.
The theorem suggests that the maximal disc of full conjugacy is strictly contained in the maximal disc of semi-conjugacy. That this can actually be the case follows from the example considered in the previous remark. 

See the work of Arrowsmith and Vivaldi \cite{Arrowsmith/Vivaldi:1994} for explicit examples of this fact, concerning power functions in the $p$-adic case.

\subsection{Proof of the theorem}

By the condition of the theorem, the conjugacy $g$ must be a power series of the form 
\[
g(x)=x +
\sum_{k=2}^{\infty}b_kx^k.
\]
By solving the equation
(\ref{semi-conjugate-equation attracting}) for formal power series we obtain
\begin{equation}\label{b_k-equation}
b_k(\lambda-\lambda^k)=\sum_{l=1}^{k-1}b_l%
(\sum\frac{l!}{\alpha_1!\cdot ...\cdot
\alpha_k!}a_1^{\alpha_1}\cdot ...\cdot a_k^{\alpha_k}),
\end{equation}
where $\alpha_1+...+\alpha_k=l$ and
$\alpha_1+2\alpha_2...+k\alpha_k=k$ and $a_1=\lambda$. Note that $|l!/(\alpha_1!\cdot ...\cdot \alpha_k!)|\leq 1$
since $l!/(\alpha_1!\cdot ...\cdot \alpha_k!)$ is an integer. Thus, by ultrametricity,
equation (\ref{b_k-equation}) yields  
\begin{equation}\label{bk-inequality}
  |b_k|\leq \frac{1}{|\lambda-\lambda^k|}\max_{l} |b_l  
  \lambda ^{\alpha_1}|\cdot| a_2^{\alpha_2}\cdot ...\cdot a_k^{\alpha_k}|,
\end{equation}
where $b_1=1$ and the maximum is taken over all solutions
$(l,\alpha_1,...,\alpha_k)$ to the system
\begin{equation}\label{index-equations}
   \left\{\begin{array}{ll}
            \alpha_1+...+\alpha_k=l,\\
            \alpha_1+2\alpha_2...+k\alpha_k=k,\\
            1\leq l\leq k-1,
        \end{array}
   \right.
\end{equation}
of equations for nonnegative integers $\alpha_i$. Note that each
$\alpha_i=\alpha_i(l,k)$ is a multi-valued function of $l$ and
$k$. For example for $k=4$ and $l=2$ we have the solution
$\alpha_1(2,4)=\alpha_3(2,4)=1$ and $\alpha_2(2,4)=0$ as well as
the solution where $\alpha_2(2,4)=2$ and
$\alpha_1(2,4)=\alpha_3(2,4)=0$.

As noted in Remark 1 above one can easily obtain the estimate 
\[
 |b_k|\leq |\lambda|^{-k}\left (1/\rho\right)^{k-1},
\]
for the coefficients of the conjugacy, using the fact that $\alpha_1 (l,k)\geq 0$. We will prove by more detailed analysis of the values of $\alpha_1(l,k)$, and induction over $k$, that in fact
\begin{equation}\label{boundbkattracting}
|b_k|\leq |\lambda|^{-\log_2 k} \left( 1/\rho\right)^{k-1} \quad \text{for } k\geq 2.
\end{equation} 
First note that by the condition of the lemma $b_1=1$. Moreover, by the condition of the lemma we also have   
\[
 |a_i| \leq \left( 1/\rho\right)^{i-1},
\]
for every $i\geq 2$. Furthermore, $a_1=\lambda$, and as $\left |\lambda\right |<1$ and $k\geq 2$ we have 
$|\lambda - \lambda^k| = |\lambda | \left |1 - \lambda^{k-1}\right |=\left |\lambda\right |$. 
Hence, $|b_2|\leq |\lambda|^{-1}\left(1/\rho\right)$.

We will proceed to prove (\ref{boundbkattracting}), using induction over $k$. 
Assume that, given $k\geq 3$ 
\begin{equation}\label{assumptionattracting}
|b_l|\leq |\lambda|^{-\log_2 l} \left( 1/\rho\right)^{l-1} \quad \text{for } l\geq 2.
\end{equation} 
Note that by (\ref{index-equations}) we have
\[
\sum_{i=2}^{k}(i-1)\alpha_i=k-l.
\]
Consequently, since $|a_i|\leq \left (1/\rho\right )^{i-1}$, we obtain
\begin{equation*}
\prod_{i=2}^k|a_i|^{\alpha_i}\leq
\prod_{i=2}^{k}\left (1/\rho\right )^{(i-1)\alpha_i}=\left (1/\rho\right)^{k-l}.
\end{equation*}
In view of the inequality (\ref{bk-inequality}) and the fact that $|\lambda -\lambda^k|=|\lambda|$ we then have
\[
 |b_k|\leq |\lambda|^{-1}\max_{1\leq l\leq k-1} |b_l  
  \lambda ^{\alpha_1(l,k)}| \left (1/\rho\right)^{k-l}.
\]
Using the assumption (\ref{assumptionattracting}) we obtain
\[
|b_k|\leq  \left (1/\rho\right)^{k-1} \max_{1\leq l\leq k-1} |\lambda| ^{\alpha_1(l,k)-1-\log_2 l}.
\]
Note the following lemma
\begin{lemma}\label{lemma-solution-2^k}
Given $k\geq 2$, let $l$ be the largest integer such that $\alpha_1(l,k)=0$. Then $l\leq k/2$.
In particular, for $n\geq 1$, we hava a solution to the equation $\alpha_1(2^n,2^{n+1})=0$.
\end{lemma}
{\noindent \bf Proof. }  Suppose that $\alpha_1(l,k)=0$.  Then, the first statement follows from the observation that by  (\ref{index-equations}) we have 
\[
k-l=\alpha_2+2\alpha_3+...+(k-1)\alpha_k\geq l.
\]
The second statement follows from the preceding and the
observation that $\alpha_2=2^n$ and
$\alpha_1=\alpha_3=...=\alpha_k=0$ is a solution of
(\ref{index-equations}) for $l=2^n$ and $k=2^{n+1}$.~$\square$

Consequently,
\[
\alpha_1(l,k)-1-\log_2 l\geq \alpha_1(l,k)-1-\log_2 (k/2)=-\log_2 k,
\]
and the assertion (\ref{boundbkattracting}) follows.

Accordingly, the radius of convergence $R_g\geq\rho$ as required. Moreover,  
$b_1=1$ and
\[
\inf_{k\geq 2}\left (1 /|b_k|\right )^{1/(k-1)}\geq \inf_{k\ge2} |\lambda|^{(\log_2k)/(k-1)}\rho
=|\lambda|\rho.
\]
It follows that the radius of injectivity $\delta_g\geq |\lambda|\rho$.
This completes the proof.

\section{Repelling fixed points}

In the repelling case we will prove the following theorem.

\begin{theorem}\label{theorem-repelling} 
Let  $K$ be a complete ultrametric field, and $f\in\mathcal{F}_{\lambda,\rho}(K)$ with 
$|\lambda|>1$. Then, there is a
unique function $g$, with
$g(0)=0$ and $g'(0)=1$, convergent and injective on the open disc $\left\{|x|<|\lambda|\gamma\right\}$,  where $\gamma=\inf_{i\geq 2}\left (|\lambda| /|a_i|\right )^{1/(i-1)}$, such that both the semi-conjugacy 
\begin{equation}\label{semi-conjugate-equation}
 g\circ f(x)=\lambda g(x),
\end{equation}
and the full conjugacy 
\[
g\circ f\circ g^{-1}(x)=\lambda x, 
\]
holds for $|x|<\delta_f$. 

\end{theorem}

\noindent
\textbf{Remark 6}.
As noted in the introduction, the existence of a convergent conjugacy in the hyperbolic case (in any dimension)   was proven by Herman and Yoccoz \cite{HermanYoccoz:1981},  generalizing results from complex dynamics.
%

\noindent
\textbf{Remark 7}.
That  $g$ converges on the disc $\left \{|x|<\gamma |\lambda | \right\}$ was also stated by Zieve 
\cite{Zieve:1996}.

\noindent
\textbf{Remark 8}.
The esitmate of the full conjugacy is maximal since $\delta_f$ is the radius of injectivity for $f$. In particular if $D_{\delta_f}(0)$ is rational. Then, $D_{\delta_f}(0)$ is the maximal disc on which the full conjugacy holds. E.g. this is the case if $f$ is a polynomial. Concerning power series, consider for example the case  $a_i=1/\lambda^{i}$ for all $i\geq 2$ given in Example 4 in section \ref{section injectivityrepelling}. Also in this case, $D_{\delta_f}(0)$ is rational since 
$\delta_f= |\lambda|$.

\begin{proof}
Again, as in the proof for attracting fixed points, we denote by $b_k$ the coefficients of the conjugacy  $g$. Recall that by the condition of the lemma
\begin{equation}\label{inf}
\gamma=\inf_{i\geq 2}\left (|\lambda| /|a_i|\right )^{1/(i-1)}.
\end{equation}
Note that if $R_f\geq \gamma$, then $\gamma =\delta_f$.  
Also note that by (\ref{inf})
\[
|a_i|\leq |\lambda |\left(1/\gamma\right)^{i-1}.
\]
We will prove that 
\[
|b_k|\leq |\lambda |^{-(k-1)}\left( 1/\gamma \right)^{k-1} \quad \text{for } k\geq 1.
\]

Again, since $g'(0)=1$ we have $b_1=1$. As in the attracting fixed point case, the coefficients of 
$g$ can then be obtained recursively according to (\ref{b_k-equation}).  As a consequence, 
$|b_2|\leq |\lambda|^{-2}|a_2|^1\leq |\lambda|^{-1}\left(1/\gamma \right)$. 
We will proceed using induction over $k$.
Given $k\geq 3$, suppose that 
\begin{equation}\label{assumption-repelling}
|b_l|\leq |\lambda |^{-(l-1)}\left( 1/\gamma \right)^{l-1} \quad \text{for } 1\leq l \leq k-1.
\end{equation}
As $|\lambda |>1$ we have $|\lambda-\lambda^k|=|\lambda|^k$. Hence, by the inequality (\ref{bk-inequality}) 
\[
|\lambda|^k  |b_k|\leq \max_{1\leq l\leq k-1} |b_l |\cdot  
|  \lambda ^{\alpha_1}\cdot a_2^{\alpha_2}\cdot ...\cdot a_k^{\alpha_k}|.
\]
Note that by (\ref{index-equations}) we have
\[
\sum_{i=2}^{k}(i-1)\alpha_i=k-l.
\]
Consequently, since $|a_i|\leq |\lambda |\left (1/\gamma\right )^{i-1}$, we obtain
\begin{equation*}
|\lambda |^{\alpha_1}\prod_{i=2}^k|a_i|^{\alpha_i}\leq
\prod_{i=1}^{k}|\lambda |^{\alpha_i}\prod_{i=2}^{k}\left (1/\gamma\right )^{(i-1)\alpha_i}
=|\lambda |^l\left (1/\gamma\right)^{k-l}.
\end{equation*}
Hence, 
\begin{equation*}
  |\lambda^{k}||b_{k} |
  \leq \max_{1\leq l\leq k-1} 
   \left [
   | b_l | |\lambda |^l\left(1/\gamma\right)^{k-l}\right ].
  \end{equation*}
The assumption (\ref{assumption-repelling}) then yields
  \[
 |\lambda^{k}||b_{k} |
  \leq \max  
   \left [
   |\lambda| \left(1/\gamma\right)^{(l-1)+(k-l)}\right ]
    =|\lambda |\left(1/\gamma\right)^{k-1}, 
\]
as required.

Hence, the radius of convergence $R_g\geq|\lambda|\gamma\geq |\lambda |\delta_f$. Recall that $f$ maps $D_{\delta_f}(0)$ one-to-one onto $D_{|\lambda|\delta_f}(0)$. Consequently, the semi-disc is of radius greater than or equal to 
$\delta_f$. Moreover,  
$b_1=1$ and
\[
\inf_{k\geq 2}\left (1 /|b_k|\right )^{1/(k-1)}\geq \inf_{k\ge2} |\lambda|\gamma
=|\lambda|\gamma.
\]
It follows that the radius of injectivity $\delta_g\geq |\lambda|\gamma\geq  |\lambda|\delta_f$.
This completes the proof.  

\end{proof}

\section*{Acknowledgement}
We would like to thank Andrei Khrennikov for fruitful discussions, and Rob Benedetto for helpful comments. In particular, we thank Juan Rivera-Letelier for corrections and 
suggestions that certainly helped us improve the presentation.



\end{document}